\documentclass[12pt,a4paper]{article}
\usepackage{graphicx}
\usepackage{amsmath,lscape}

\usepackage{multirow}

\title{Average treatment effect estimation via random recursive partitioning}
\author{Giuseppe Porro\footnote{{\it corresponding author}: Department of Economics and Statistics, University of Trieste, P.le Europa 1, I-34127 Trieste, Italy. E-mail:giuseppe.porro@econ.units.it}\\ Stefano Maria Iacus\footnote{Department of Economics, University of Milan, Via Conservatorio 7, I-20122 Milano, Italy. E-mail: stefano.iacus@unimi.it}}
 \date{September 2004}
\begin{document}

\maketitle

\begin{abstract}

A new matching method is proposed for the estimation of the average
treatment effect of social policy interventions (e.g., training programs or
health care measures). Given an outcome variable, a treatment and a set of pre-treatment covariates,
the method is based  on the examination of random recursive partitions of the space 
of covariates using regression trees. A regression tree is
grown either on the treated  or on the untreated individuals {\it only}
using as response variable a  random permutation  of the indexes  1\ldots$n$ ($n$ being the number of units involved), while the indexes for the other group are predicted using this tree. The procedure is replicated in order to rule out the effect of specific permutations. The average treatment effect is estimated in each tree by matching treated and untreated in the same terminal nodes. The final estimator of the average treatment effect is obtained by averaging on all the trees grown.
The method does not require any specific model assumption apart from the tree's complexity, which does not affect the  estimator though.
We show that this method is either an instrument to check whether  two samples can be matched (by any method) and, when this is feasible, to obtain reliable estimates of the average treatment effect. We further propose a graphical tool to inspect the quality of the match.
The method has been applied to the National Supported Work Demonstration data, previously analyzed by Lalonde (1986) and others.

\end{abstract}


\noindent
{\bf J.E.L. Classification:} C14;\ C40;\ C81

\newpage

\section{Introduction}

A wide category of estimators has been developed in the last decade to
evaluate the effects of medical, epidemiological and social policy interventions (for instance, a
training program) over the individuals (for a survey on methods and applications see Rosenbaum, 1995 and Rubin, 2003). Matching
estimators represent a relevant class in this category:
these estimators aim to combine ({\it match}) individuals who have been subject to the
intervention (so forth, {\it treated}) and individuals with similar pre-treatment characteristics who have not
been exposed to it ({\it untreated} or {\it controls}), in order to estimate the
effect of the external intervention as the difference in the value of an outcome
variable.

From a technical viewpoint, a remarkable obstacle in obtaining a satisfying match is constituted by the ``curse of dimensionality'', i.e. by the dimension of the covariate set.
This set  should
be as large as possible, in order to include all the relevant informations
about the individuals (Heckman {\it et al.} 1997, 1998), but increasing the number of covariates makes the match a complicated task. 

Two main approaches have been developed to solve the matching problem through a one-dimensional measurement:  the propensity score matching method (PSM) and matching method based on distances (DM).
The PSM method makes use of the notion of {\it propensity score}, which is defined as the probability to be exposed to the treatment, conditional on the covariates. Treated and controls are matched on the basis of their scores according to different criteria: stratification, nearest neighbor, radius, kernel (see Smith and Todd, 2004a, for a review).  The DM method makes use of  specific distances (e.g. Mahalanobis)  to match treated and untreated (see e.g. Rubin, 1980, Abadie and Imbens, 2004).

In this paper we propose a matching method based on the exploration of random partitions of the space of covariates. Treated and controls are considered similar, i.e.  matched,  when belonging to the same subset in a partition. As it will be clear in the following, our  technique is unaffected by the dimensionality problem. 

More precisely, we make use of a non standard application of regression trees. The CART (classification and regression trees) methodology
has been introduced in the literature on matching estimators as one of the
alternatives to parametric models in the assignment of propensity scores (see e.g. Rubin, 2003 and Ho {\it et al.}, 2004) or to directly match individuals inside the tree (Stone {\it et al.}, 1995). 
In our approach the use of regression trees is indeed different: our
technique exploits the ability of trees to partition the space of
covariates.

The methodology proposed here is to grow a regression tree only on one group, for example the treated, in the following way: we assign each unit a progressive number (label) and use it as a response variable.  We grow the tree till it has only one (or at most few) units in the terminal nodes (leaves). Then we use the tree to predict the labels for the units of the other group  (in this example, the controls). 
Units (treated and controls) are matched if they are in the same leaf and only if the balancing property for their covariates is met.
As the resulting tree (i.e. the partition of the covariate space) depends strongly on the initial assignment of the labels  (see Section \ref{sec:cart}), we operate several permutations of these labels  and then take the average of all the results to get the final treatment effect estimate.

To shed light on the properties of the method, we use the National Supported Work Demonstration (NSW) data coming from a training program originally analyzed by Lalonde (1986) but also by a large number of studies aiming to test the performance of different methodologies of evaluation (Dehejia and Wahba, 1999, 2002; Becker and Ichino, 2002; Smith and Tood, 2004a, 2004b; Deheja 2004; Abadie and Imbens, 2004)

This dataset is peculiar in its previous analyses  pointed out that ``the question is not which estimator is the best estimator always and everywhere'' (Smith and Todd, 2004a);
the task of the investigation should be to provide a tool able to signal when the matching methods can be successful or, on the contrary, alternative methodologies ought to be applied. The procedure proposed in this paper is a possible candidate to perform the task; moreover, when matching can be applied, the estimators obtained with this procedure  are either normally distributed, robust with respect to the complexity of the tree and capable to reduce the bias.

The paper is organized as follows: in Section \ref{sec:match} we introduce the notation and define the object of estimation. Section \ref{sec:cart} is devoted to a brief description of the CART methodology which is functional to Section \ref{sec:idea}, where  the idea behind the proposed technique is illustrated. In Sections \ref{sec:res} and \ref{sec:algo} we present in details the algorithmic flow of the procedure  and our empirical results to show the above mentioned properties of the method. Section \ref{sec:proxy}
introduces a graphical tool which is useful to asses the quality of the match. Most of the tables and all the figures can be found at the end of the paper.

\section{The matching problem}\label{sec:match}
Matching estimators have been widely used in the estimation of the average treatment effect (ATE) of a binary treatment on a continuous scalar outcome. For individual $i = 1, \dots, N$, let $(Y_i^T, Y_i^C)$ denote the two potential outcomes,  $Y_i^C$ being the outcome of individual $i$ when he is not exposed to the treatment and $Y_i^T$ the outcome of individual $i$ when he is exposed to the treatment to estimate $Y_i^C$. For instance, the treatment may be participation in a job training program and the outcome may be the wage. If both $Y_i^C$ and $Y_i^T$ were observable, then the effect of the treatment on $i$ would be simply $Y_i^T - Y_i^C$. The root of the problem is that only one of the two outcomes is observed, whilst the counterfactual is to be estimated.

Generally, the object of interest in applications is the average treatment effect on the subpopulation of the $N_T$ treated subjects (ATT). Let $\tau$ be the ATT, then $\tau$ can be written as
$$\tau=\frac{1}{N_T} \sum_i (Y_i^T - Y_i^C)$$ 
As said, the first problem  in practice is to estimate the unobserved outcome, $Y_i^C$ for individual $i$  who was exposed to the treatment.  If the assignment to the treatment were random, then one would use the average outcome of some similar individuals who were not exposed to the treatment. This is the basic idea behind matching estimators. For each $i$, matching estimators impute to the missing outcome the average outcome of untreated individuals similar, up to the covariates, to the treated ones.

To ensure that the matching estimators identify and consistently estimate the treatment effects of interest, we assume that: a)  (\textit{unconfoundedness}) assignment to treatment is independent of the outcomes, conditional on the covariates; b) (\textit{overlap}) the probability of assignment is bounded away from zero and one (see Rosenbaum and Rubin, 1983; Abadie and Imbens, 2004, and references there in).

As mentioned in the Introduction, the problem with matching is that the number of covariates and their nature makes usually hard to provide an exact match between treated and control units.

\section{About classification and regression trees}\label{sec:cart}
 Classification and Regression Trees (CART) have been proposed by Breiman {\sl et al.} (1984) as a classifier or predictor of a response variable (either qualitative or quantitative), conditionally on a set of observed covariates,  which collects  sample units in groups  as much homogeneous as possible with respect to the response variable.  
 The main assumption in CART is that the set of covariates $\mathcal X$ admits a partition and the tree is just a representation of this partition. To our end this means that the space $\mathcal X$ is divided into cells or strata where we can match treated and control units.

Ingredients for growing a tree are the space of covariates $\mathcal X$, the response variable $Y$ and a homogeneity criterion (e.g. the deviance or the Gini index).  One of the most commonly used methods to grow a tree is the ``one-step lookahead'' tree construction with binary splits (see Clark and Pregibon, 1992). 
Given the set $\mathcal X$ we say that  its partition can be represented as a tree or, more precisely, by its terminal nodes called \textsl{leaves}. Data are then subdivided in groups according to the partition of $\mathcal X$.

In the first step of the procedure the covariate space $\mathcal X$ is not partitioned and all the data are included in the \textit{root} node. The root is then splitted  with respect to one covariate $X_j$ into two subsets such that  $X_j > x$ and $X_j \leq x$ (in case $X_j$ is a continuous variable, but similar methods are conceived for discrete, qualitative and ordered variables). The variable $X_j$ is chosen among  all the $k$ covariates $X_1, X_2, \ldots, X_k$ in such a way that the reduction in deviance inside each node is the maximum achievable. This procedure is iterated for each newly created node.
The tree construction is stopped when a minimum number of observations per leaf is reached or when the additional reduction in deviance is considered too small.

Any new observation can be classified (or its value predicted)  by dropping it down through the tree. Note that even observations with missing values for some covariates can be classified this way.

The use of CART to replace parametric models (probit or logit) to estimate the propensity score is not unknown to the literature (see e.g.  Ho {\it et al.}, 2004). 
CART is also used to  directly  match treated and controls units inside the leaves (see e.g. Stone {\it et al.}, 1995, among the others).
These approaches are used to solve the problem of model specification that is always difficult to justify in practice (see Dehejia, 2004, for an account on the sensitivity of the estimate due to model specification). 

CART  technique is also seen as  a variable selection method and it is rather useful compared to parametric modelling because interaction between variables and polynomial transforms are handled automatically.

\section{Random partition of the space of covariates using CART}\label{sec:idea}
In this section we present a new approach to the use of CART. As previously noticed, CART methodology generates  partitions of  the space $\mathcal X$  maximizing homogeneity inside the leaves with respect to the response variable.
Our proposal  is to build a tree {\it only} on the treated and to grow it up to a level of complexity sufficiently high, such that  each terminal leaf contains at most few treated. 
This produces a partition of $\mathcal X$ that reflects the stratification structure of the treated. 
To this end,  being $n_T$ the sample size of the treated, 
we  assign to the treated a response variable which is a sequence of numbers  from 1 to $n_T$.

This tree is used to assign controls to leaves. Treated and controls belonging to the same leaf are then directly matched.  Formally the ATT is  estimated as follows:
\begin{equation}
\begin{aligned}
\widehat{att} &=\frac{1}{n_T} \sum_{i\in T}\left(Y_i^T - \frac{1}{|C_i|}  \sum_{j\in C_i}  w_{ij} Y_j^C \right)\\
&= \frac{1}{n_T} \sum_{i\in T} Y_i^T - \frac{1}{n_T}   \sum_{i\in T} \sum_{j\in C_i}  \frac{w_{ij}}{|C_i|} Y_j^C
\label{eq1}
\end{aligned}
\end{equation}
where  $w_{ij}=1$ if treated $i$ and control $j$ are matched in a leaf, otherwise $w_{ij}=0$.  $T$ is the set of the indices for the treated, $C_i = \{ j : w_{ij}=1\}$, $i\in T$ and $|C_i|$ is the number of elements in set $C_i$.
The variance of $\widehat{att}$ can be directly calculated:
\begin{equation}
{\rm Var} (\widehat{att}) = \left(\frac{1}{n_T}\right)^2 \sum_{i\in T}\left( {\rm Var}(Y_i^T) + \sum_{j\in C}\left(\frac{w_{ij}}{|C_i|}\right)^2 
 {\rm Var}(Y_i^C) \right)
\label{eq2}
\end{equation}
being $C$ the set of indices for the controls. Further, pose  $w_{ij}/|C_i|=0$ if  $w_{ij}=0$ by definition.
Provided we are given a consistent estimator $\hat\mu_0(X)$ of $\mu_0(X) = E(Y^C|X)$, we can also adjust the bias for the difference-in-covariates (see Abadie and Imbens, 2004) obtaining the following adjusted version of
 equation \eqref{eq1}:
\begin{equation}
\widehat{att}'=
\frac{1}{n_T} \sum_{i\in T}\left((Y_i^T -\hat\mu_0(X_i)) - \frac{1}{|C_i|}  \sum_{j\in C_i}  w_{ij} (Y_j^C -\hat\mu_0(X_j)) \right)
\label{eq3}
\end{equation}
and the variance of $\widehat{att}'$ is the following:
\begin{equation}
\begin{aligned}
{\rm Var} (\widehat{att}') &= \left(\frac{1}{n_T}\right)^2 \sum_{i\in T}\left( {\rm Var}(Y_i^T) + \sum_{j\in C}\left(\frac{w_{ij}}{|C_i|}\right)^2 
 {\rm Var}(Y_i^C) \right)\\
 &+
 \left(\frac{1}{n_T}\right)^2 \sum_{i\in T}\left( {\rm Var}(\hat\mu_0(X_i)) + \sum_{j\in C}\left(\frac{w_{ij}}{|C_i|}\right)^2 
 {\rm Var}(\hat\mu_0(X_j)) \right)
\end{aligned}
\label{eq4}
\end{equation}
\subsection{Generating random partitions}
The tree structure is strongly dependent on the assignment of the values ($1$, \ldots, $n_T$) of the response variable to the treated; so does the ATT estimate.
In order to marginalize this dependence,  we replicate the tree construction randomly permuting the set  ($1$, \ldots, $n_T$) of the values of the response variable assigned to each treated unit\footnote{In other direct matching estimators, the same argument on order dependency applies as the order chosen to match the individuals affects the final estimator of the ATT (see e.g. Abadie and Imbens, 2004, D'Agostino, 1998).}.
The final estimator of the ATT will be the average of the ATT's obtained in each replication.
The number of possible permutations is $n_T!$ which is also the maximum number of the significant partitions of the space $\mathcal X$.

At each replication, the balancing property on the covariates has to be tested inside each leaf. When this property is not met, the treated and control units involved are excluded from the matching.

Moreover, to increase the number of matched treated in each replication, we generate a subsequence of trees to match the residual treated: first a tree is grown and matching is taken over. If not all the treated have been matched, we keep track of these and we grow another tree using a different permutation of ($1, \ldots, n_T$). At this second step we now match the remaining treated only. The procedure is iterated until either all the treated units have been matched or  a prescribed maximum number of iterations has been reached.This can be done with or without replacement for the controls.
So we have $R$ replications and for each replication we estimate the ATT as in \eqref{eq1} (or \eqref{eq3}). Define $\hat\tau^{(k)}$ the ATT estimator at $k$-th replication ($k=1, \ldots, R$), as follows:
\begin{equation}
\hat\tau^{(k)} =\frac{1}{n_T} \sum_{i\in T}\left(Y_i^T - \frac{1}{\left|C_i^{(k)}\right|}  \sum_{j\in C_i^{(k)}}  w_{ij}^{(k)} Y_j^C \right)
\label{eq1a}
\end{equation}
The final ATT estimator is defined as the average over all the replications:
$$
\hat\tau = \frac{1}{R}\sum_{k=1}^R \hat\tau^{(k)}
$$
At each replication, the set of weights $w_{ij}^{(k)}$ can be represented as  a matrix $W^{(k)}$ of matches.
Define the {\it proximity matrix} ${\bf P}$  as follows:
\begin{equation}
{\bf P} = \sum_{k=1}^R W^{(k)} = \left[p_{ij} = \sum_{k=1}^R w_{ij}^{(k)}\right]
\label{eq:proxy}
\end{equation}
This matrix, examined in detail  in Section \ref{sec:proxy}, contains information on 
the quality of the match between treated and controls. In fact, each row $i$ of $\bf P$  tells how many times, over $R$ replications, a treated unit $i$ has been matched with each control. Even if this is not a \textsl{distance} matrix, it can be used as a starting point for calculating DM-estimators.

\subsection{Further improvements}\label{sec:switched}
As usual in the applications, $n_T << n_C$. So it might happen that the tree grown up on the treated contains a high number of controls in each terminal leaf. This may cause the balancing property to fail in a relevant number of cases even in spite of a long subsequence of iterations. This implies that a small number of treated are included in the matching and generates additional bias in the estimates.

If this is the case, we propose to grow a tree \textit{only} on the controls. One might expect to find a greater number of terminal leaves with a few controls and treated per leaf by construction.
In such a tree, the balancing property should be met more frequently.

On balanced samples both kind of approaches generate similar partitions. On the contrary, for unbalanced samples, this alternative procedure may be effective.

\section{Empirical results}\label{sec:res}
In this section we present an application of the above procedure analyzing, once again, the well know benchmark example of the NSW data from Lalonde (1986). In the view of reproducible research, all the examples, including data sets, software and scripts, are available at the web page {\tt http://www.economia.unimi.it/rtree}. The software we use is the open source statistical environment called \texttt{R} which has recently increased its popularity also among econometricians (see for example, among the others,  Kings' projects \texttt{Zelig} and  \texttt{MatchIt}, also including specific matching routines; Kosuke {\it et al.}, 2004;  Ho {\it et al.}, 2004).

\subsection{The data}
The National Supported Work (NSW) data comes from LalondeÕs (1986) seminal study on the comparison between experimental and non-experimental methods for the evaluation of causal effects. The data contain the results of a training program: the outcome of interest is the real earnings of workers in 1978, the treatment  is the participation to the program. Control variables are age, years of education, two dummies for ethnic groups: black and hispanic, a dummy for the marital status, one dummy to register the posses of a high school degree, earnings in 1975 and earnings in 1974. The set contains 297 individuals exposed to the treatment and 425 control units. This is an experimental sample, named LL through the text. In their 1999 and 2002 papers, Dehejia and Wahba select a subset of LL on the basis of 1974 earnings, restricting the sample to  185 treated and 260 controls. We call this sample DW. Further, Smith and Todd (2004a) suggest a more consistent selection on the basis of 1974 earnings of individuals from the LL sample, limiting the number of treated to 108 and to 142 for the controls. We refer to this set a ST.
Along with these experimental samples, it is common practice to use a non experimental sample of 2490 controls coming from the  Population Survey of Income Dynamics. These data,  called PSID in the following, are used generally to prove the ability of PSM or DM methods in matching experimental and non-experimental data (see e.g. Dehejia and Wahba, 1999, 2002; Abadie and Imbens, 2004) or to show, on the contrary, that LL (and its sub-samples) are not comparable with PSID data (Smith and Todd, 2004a).
In the following we call ``naive target'' the difference in mean of the outcome variable (earnings in 1978) of treated and controls in the experimental samples (LL, DW and ST). This is considered to be the benchmark value for the ATT. 

\subsection{The implementation of the random tree approach}
In our approach, trees are grown to their maximum complexity using the recursive partition method as suggested by Breiman {\it et al.} (1984). The only varying parameter  is the minimum split (just {\it split} in the tables), which is the minimum number of observations that must be included in a leaf for a split to be attempted (this means that with a split equal to 2, in absence of ties, the resulting tree has one observation per leaf). Like in Abadie and Imbens (2004), there is only one parameter the user needs to choose during the analysis and yet the estimates do not show particular sensitivity with respect to it.
We made 250 replications, i.e. we search for 250 different partitions of space of the covariates, allowing for a maximum number of 50 iterations per replications (which means a maximum of 250*50 =12500 trees grown; for the actual number of trees grown see Table \ref{tab4}). A small number of iterations per replication means that matching is easily attainable.
\begin{table}[h]
\par
\hrule
\begin{center}
\begin{tabular}{r|cccc||cccc}
\multicolumn{1}{c}{}&\multicolumn{4}{c}{Trees grown on the treated} &\multicolumn{4}{c}{Trees grown on the controls} \\
split & LL & DW & ST & DWvsPSID & LL & DW & ST & DWvsPSID\\
\hline
50 & 5 & 4 &2 &43&6&4&2&42\\
32 & 6& 5 & 4&49&7&5&5&39\\
20 & 10& 6 &5&50&8&7&6&35\\
16 &15 &12&9&50&10&8&7&36\\
8&47&45&34&50&15&14&11&37\\
4&50&50&50&50&17&17&8&41\\
2&50&50&50&50&45&43&37&50\\
\end{tabular}
\end{center}
\caption{\small Average number of iterations per replication. In every sample studied, the average number of trees grown varies from 2*250=500 to 50*250=12500, being 250 the number of replications. The higher the number of iterations the harder  the matching. For matchable samples the number of iterations increases with the split.}
\label{tab4}
\hrule
\par
\end{table}

Contrary to some DM and PSM implementations in the literature, we still check for the balancing property in each leaf of the trees using a $t$ test for quantitative covariates
and a $\chi^2$ test for qualitative or dicotomus variables using a significance level of 0.005.
The difference-in-covariates correction is applied using a standard regression tree to estimate $\mu_0(X)$ (see Section \ref{sec:idea}). 

\subsection{Results on experimental data}\label{sec:exp}
In Tables \ref{tab1} and \ref{tab2} we present the results for the experimental samples LL, DW and ST.
From the empirical analysis in Table \ref{tab1} it emerges that ATT estimators $\hat\tau$ and $\hat\tau'$ are both close to the target and stable (i.e. non sensitive to the split parameter) in all the samples. The estimated standard deviation of the population in equation \eqref{eq2} is denoted by $\hat\sigma_{\tau}$ (and $\hat\sigma_{\tau'}$ respectively) and is also close to the target (in this case the pooled standard deviation).
Table \ref{tab1} and \ref{tab2} also report the following informations: the mean number of treated and controls matched per replication, the number of trees in 250 replications that matched more than 95\% of the treated units in the sample (``\% o.t.'' in the tables, namely ``over the 95\% threshold''), the 95\% Monte Carlo confidence intervals for $\tau$ and $\tau'$ and the results SW and SW' for the Shapiro-Wilk normality test (see e.g. Royston, 1995) on the distribution of $\tau$ and $\tau'$: a significant value for the test (one or more bullets in the tables) means non normality.

The difference between Tables \ref{tab1} and \ref{tab2} is in the way the trees are grown. Usually ``switched'' trees (trees grown on the controls only) of Table \ref{tab2} match almost all the treated even at very small split. This is due to the fact that the control group is numerically bigger than the treated one, so that the switched tree has a number of leaves higher than the  number of treated units and match frequently occurs one-to-one. Conversely, on ``straight'' trees (trees grown on the treated only) controls tend to group on leaves where only one treated is present, often causing the violation of the balancing property. 
In Table \ref{tab1} one can see that when the percentage of trees that match at least 95\% of treated is low, the estimators do not exhibit a good performance. So this percentage can be assumed as a quality indicator of the match. In principle, the evaluation of $\tau$ requires to match {\it all} the treated units: a poor match, i.e. the systematical exclusion of some treated individuals from the matching process, forces the ATT estimator to neglect a part of the treatment effect and hence introduces a bias in the estimates.
Notice that from this point of view our procedure is rather selective compared to DM methods or to some implementations of the PSM methods, because we require  the balancing property to be met  inside each leaf.
Also note that the confidence intervals for $\tau'$ are usually a bit smaller than the analogous for $\tau$, meaning that -- despite the population variance is affected by the difference-in-covariates correction --this is not likely to happen for the standard error of the estimator.
Remarkably, the results highly agree in both tables.

\subsection{Results on non-experimental data}\label{sec:nonexp}

Table \ref{tab3} shows the results obtained by matching DW treated units and PSID
controls. The use of non-experimental controls to evaluate the average
effect of treatment of NSW treated workers has been attempted in many
studies and the conclusion of Smith and Todd (2004a) is that low bias
estimates cannot be achieved by matching estimators when using samples
coming from different contexts (different labour markets, in  this case)
and relying on non-homogeneous measures of the outcome variable. Therefore,
they conclude that NSW data  and PSID cannot reach a good
match. The same evidence arise in Abadie and Imbens (2004) where the authors try to match the same two samples. We obtain results qualitatively similar to these studies.

As in Smith and Todd (2004a) and Abadie and Imbens (2004), the results in Table \ref{tab3} show a bad quality of the
matching and the consequent unreliability of the ATT estimates. Differently to the quoted literature, the method proposed here gives a clear view that something is really wrong in matching these two samples (even in absence of a target).

At first, notice that the trees that match more than 95\% of the
treated (``\% o.t.'') never exceed 59\% of the total and that the average
number of treated units matched in each tree is never higher than 177 (over
185 treated in DW sample). As previously discussed, having a higher number
of leaves, the ``switched'' trees achieve better values of the ``\% o.t.''\
indicator and match, on average, a higher number of treated per tree,
together with a lower average number of controls. The average number of
iterations - never lower than 35 and often equal to the maximum - confirms
the difficulty of matching the two samples (see Table \ref{tab4}). Moreover, the
number of iterations appears to be higher the lower is the split parameter,
i.e. the lower is the ``\% o.t.''\ indicator.

On the whole, all the indicators seem concordant in pointing out that DW
treated and PSID controls cannot achieve a good quality match, and this
should suggest that these two samples are not to be used to evaluate the average treatment effect through a matching estimator. Examining the values
of $\hat\tau$\ and $\hat\tau'$ in Table \ref{tab3} one can in fact observe that the estimated ATT is quite far
from the naive target (1794), particularly before the
difference-in-covariates correction.

It is more considerable, however, that many features of the estimates would
confirm this unreliability even in the absence of a benchmark. The main
evidence is the remarkable difference between the ATT estimate before and
after the difference-in-covariates adjustment: in our case the signs of $\hat\tau$ and $\hat\tau'$
are sistematically different. The difference-in-covariates adjustment
brings the estimates closer to the target, especially when the ATT is
evaluated from the ``switched'' trees: this indicates that, due to the
difficulty in matching the two samples, the lack of balance in the
covariates inside the leaves is still considerable when $\tau$ is estimated. As one may expect (see \S \ref{sec:switched}), when the two
samples cannot be easily matched, the ``switched'' trees provide lower bias
estimates. The estimated variance of the ATT is particularly high, both with
trees grown on the treated and with trees grown on the controls. Also the
variance of the estimator assumes high values, so that the confidence
intervals we obtain include the target value of the ATT, at least after the
difference-in-covariates adjustment (and with a partial exception when 
split=2), only due to the wideness of the intervals themselves. Lastly, the
Shapiro-Wilk test shows that the estimators seldom have a normal
distribution.

\section{Summary of the method}\label{sec:algo}
In the previous section we  showed the performance of the proposed tool in two typical situations of  average treatment effect estimation on experimental (\S \ref{sec:exp}) and non-experimental data (\S \ref{sec:nonexp}). The results seem to be able to answer to the main questions: ``are my data really comparable in terms of the matching procedures?''; if so, ``can I reliably estimate the average treatment effect?'' More importantly, both questions can be addressed even in absence of a known target and without needing to introduce model assumptions hardly justifiable in practical situations.
We would like to summarize here the properties of the method and give an outline on how to use this tool in applied analyses.

\subsection*{The characteristics of the method}
\begin{itemize}
\item  the method, as desirable for any method, does not require to know the target in advance to provide conclusions/evidence on the data;
\item the method does not need to specify  any additional model assumption;
\item the method is not sensitive to the only parameter (the split) the user can specify (this is actually the only, so to say,  ``assumption'' we put on the complexity of the grown trees, i.e. on the model);
\item explicit indications of low quality of the match are given by the method, in terms of difficulty to match treated and controls (see Table \ref{tab4}) and number of trees that match almost all (at least 95\% of) the treated (see Tables \ref{tab1}, \ref{tab2}, \ref{tab3});
\item  the estimating procedure provides information on the reliability of the
estimates (i.e. when treated and controls are correctly matched the difference-in-covariate correction should not considerably affect the ATT estimate);
\item the procedure is available to be used as an open source free software and can be easily implemented for other software environments like Stata, Matlab etc.
\item the method has at least one limitation: no formal theory is available yet. Conversely, one can think of our approach as a non-parametric stratified permutation test approach to the problem of average treatment effect estimation. This seems to be a promising way --  still under investigation -- for proving formal properties of the ATT estimator.
\item even if not stressed in the text, CART is a method known to be robust to missing values in the dataset and this can be a important advantage in some non-experimental situations.
\end{itemize}
\subsection*{How to use the procedure}
We conclude this section illustrating  a sort of step-by-step guide, in algorithmic form, for using the method in applied research.
\begin{enumerate}
\item[]{\it Step 0.} initialization: choose {\tt k}, the number of replications, to be at least {\tt 100}, and set the flag {\tt iter} to {\tt FALSE};
\item[]{\it Step 1.} run {\tt k} replications using a maximum number of 50 iterations on the straight and switched trees. Let the split vary from, say, 50 to 2.
\item[]{\it Step 2.} observe the following indexes for the different values of the split parameter: a) the number of partitions/replications that match at least 95\% of the treated; b) the average number of iterations per replication; c) the values of $\hat\tau$ and $\hat\tau'$; 
\item[]{\it Step 3.} Three cases may happen:
\begin{itemize}
\item[]{\it Case i:} a) is low, the match is not successfully realized. In this case you should notice high values in b) and a substantial effect on the difference-in-covariates adjustment in c). Match is not the tool to analyze this dataset.  {\bf Stop}.
\item[]{\it Case ii:} a) is high, b) is low and c) present stable results then c) gives you reliable informations on the true value of the ATT. 
\begin{itemize}
\item[]{$\bullet$} If {\tt iter=FALSE} set it to {\tt TRUE}. Go to {\it Step 1} using a higher number  {\tt k} of replications to obtain more precise ATT estimates (for example in terms of confidence intervals). You don't need to grow both straight and switched trees at this point: just build the trees on the subset with the higher number of observations. You can also drop values of the split (usually the low values) corresponding to lower values of a). 
\item[]{$\bullet$} If {\tt iter=TRUE} then  you can rely on the ATT estimates. {\bf Stop}
\end{itemize}
\item[]{\it Case iii:} you cannot draw sharp conclusions observing the values in a), b) and c) for different values of the split. 
\begin{itemize}
\item[]{$\bullet$} If {\tt iter=FALSE} set it to {\tt TRUE}, go to {\it Step 1} using a higher number {\tt k} of replications, possibly operating a selection on the split values.
\item[]{$\bullet$}  If {\tt iter=TRUE} then you have no reliable ATT estimates: matching is not the tool to analyze this dataset. {\bf Stop}. 
\end{itemize}
\end{itemize}
\end{enumerate}
Given the random nature of the method, one iteration ({\tt iter=TRUE}) of the protocol with an increased number {\tt k} of replications is always recommendable to eliminate the influence of the initial seed of the random generator used. As already mentioned in the text, with {\tt k=250} replications and a maximum number of iterations per tree equal to 50, the total number of trees generated (i.e. the total number of partitions explored) results to be appreciably high, up to 12500.

\section{The proximity matrix and a related simple ATT estimator}\label{sec:proxy}
\begin{figure}
\begin{center}
\begin{tabular}{c}
\includegraphics[height=7cm]{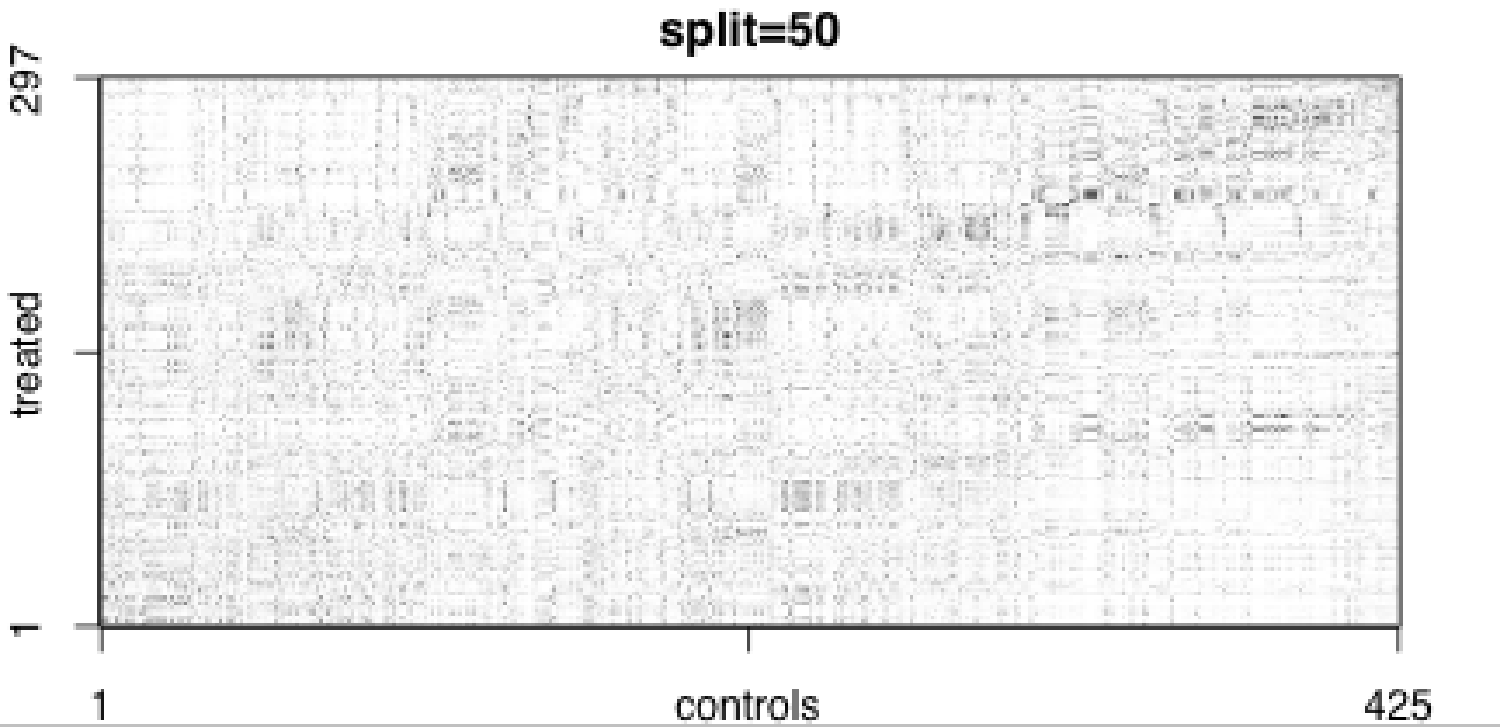}\\
\includegraphics[height=7cm]{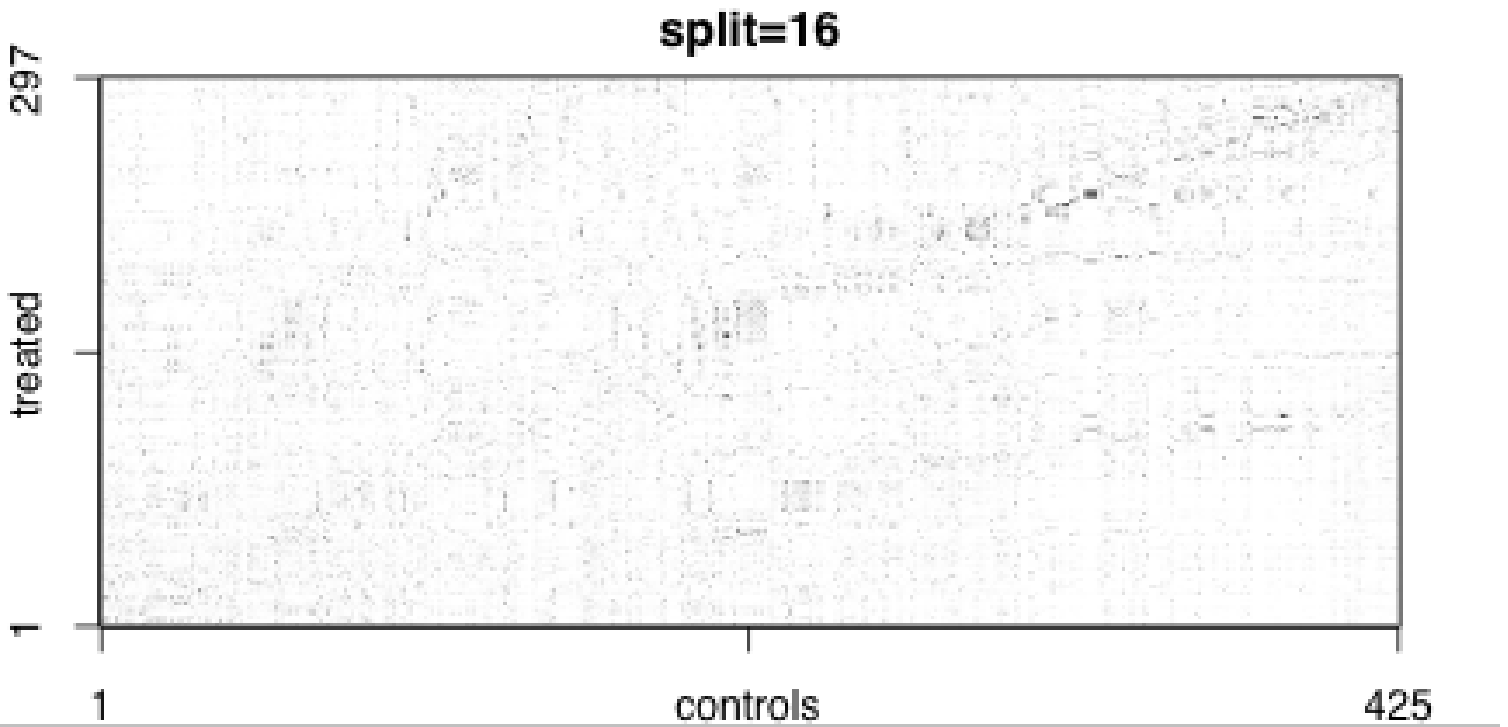}\\
\includegraphics[height=7cm]{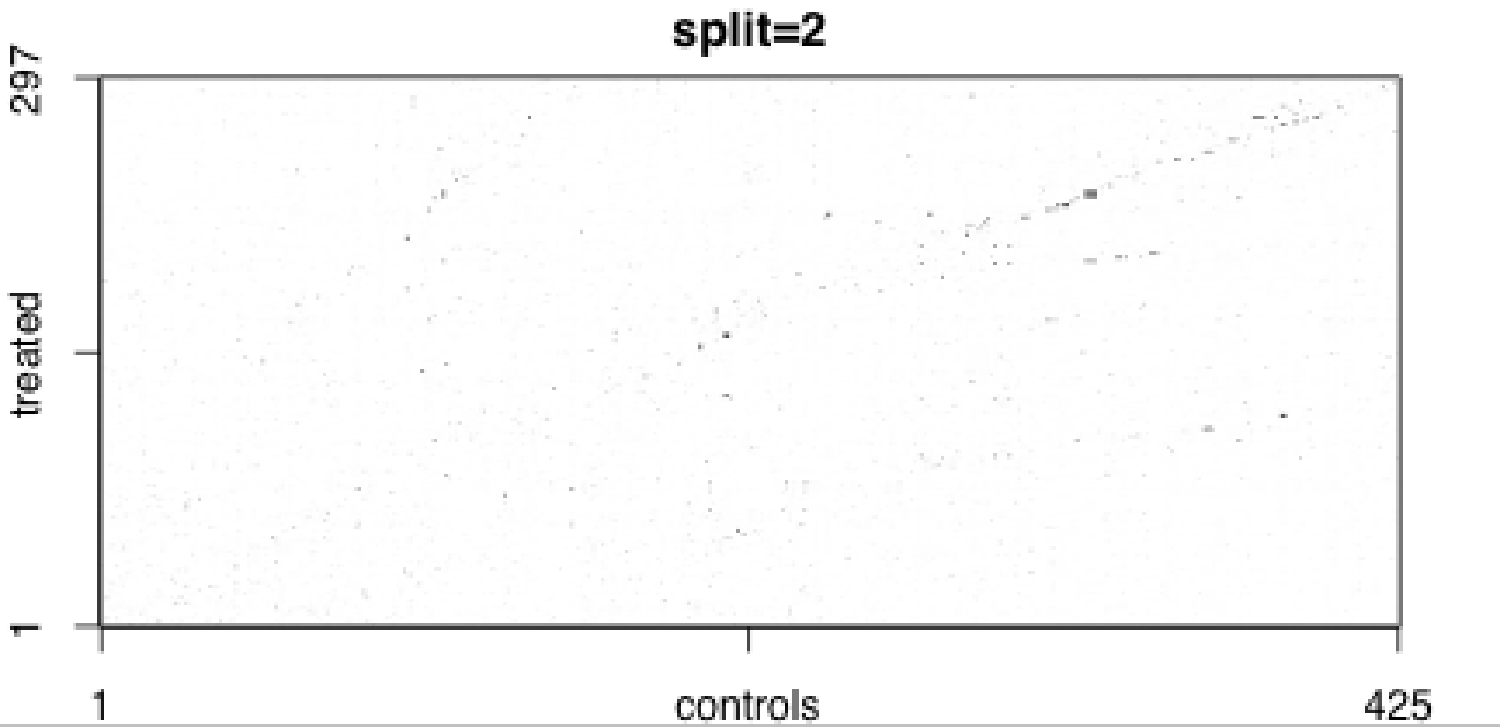}
\end{tabular}
\end{center}
\caption{\small Proximity matrix for random partitions built on the controls. LL experimental dataset for different split values. Large values in the split parameter produce spurious matches that tend to vanish as the split values decreases. In the last matrix, only ``true'' matches survive.}
\label{fig1}
\end{figure}

The proximity matrix $\bf P$ defined in \eqref{eq:proxy} summarizes the number of matches realized along the $R$ replications. The scalar $p_{ij}$ of $\bf P$ reports the number of matches involving treated $i$ and control $j$.
The graphical representation of $\bf P$ allows for an easy inspection of the quality of the match. Figure \ref{fig1}  is a  representation of matrix $\bf P$ for split values 2, 16 and 50 for the LL experimental sample. Treated units are represented on the $y$ axis while the control units are on the $x$ axis. The intensity of the spots is proportional to the corresponding $p_{ij}$: the darker the spot, the higher the number of times the match between treated $i$ and control $j$ is realized. 
Given a split=2, paired samples are expected to return an image with one spot per line. As long as a treated units matches more than a single control, the corresponding line reveals several points.
This  occurs either when imposing  larger values of the split parameter (spurious matches) or when using unpaired samples (a treated has more than one truly similar individual among the controls), or for both reasons.
Therefore, the two conditions for a good match are: a) having at least one spot per line and b) having spots as dark as possible.
Images with many faint spots  and few dark spots  are the result of a low quality 
 match realization.
When samples can be successfully matched, reducing the split value should generate  a  cleaner and sharper image, as spurious matches are removed when the split decreases and the darkness of the spots representing true matches is not affected.

\begin{figure}
\begin{center}
\begin{tabular}{c}
\includegraphics[height=7cm]{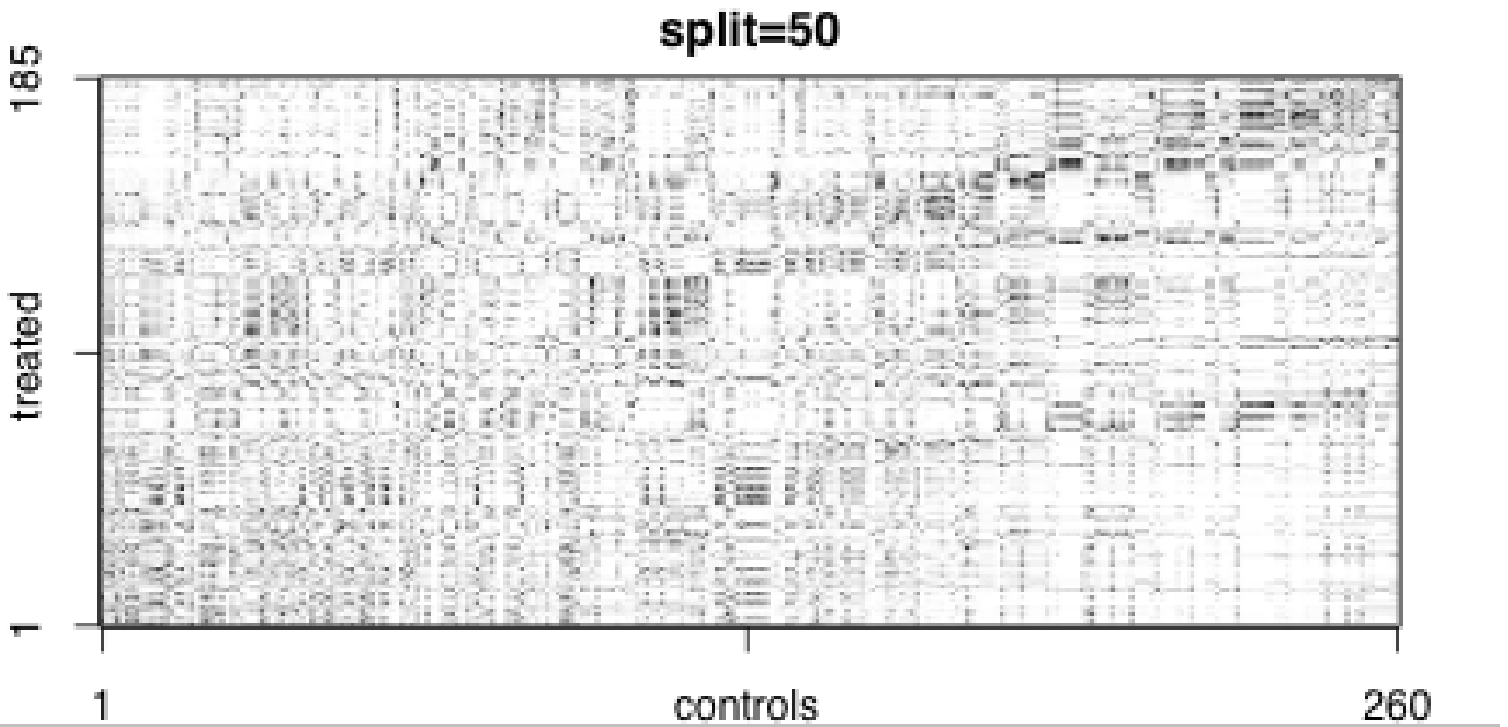}\\
\includegraphics[height=7cm]{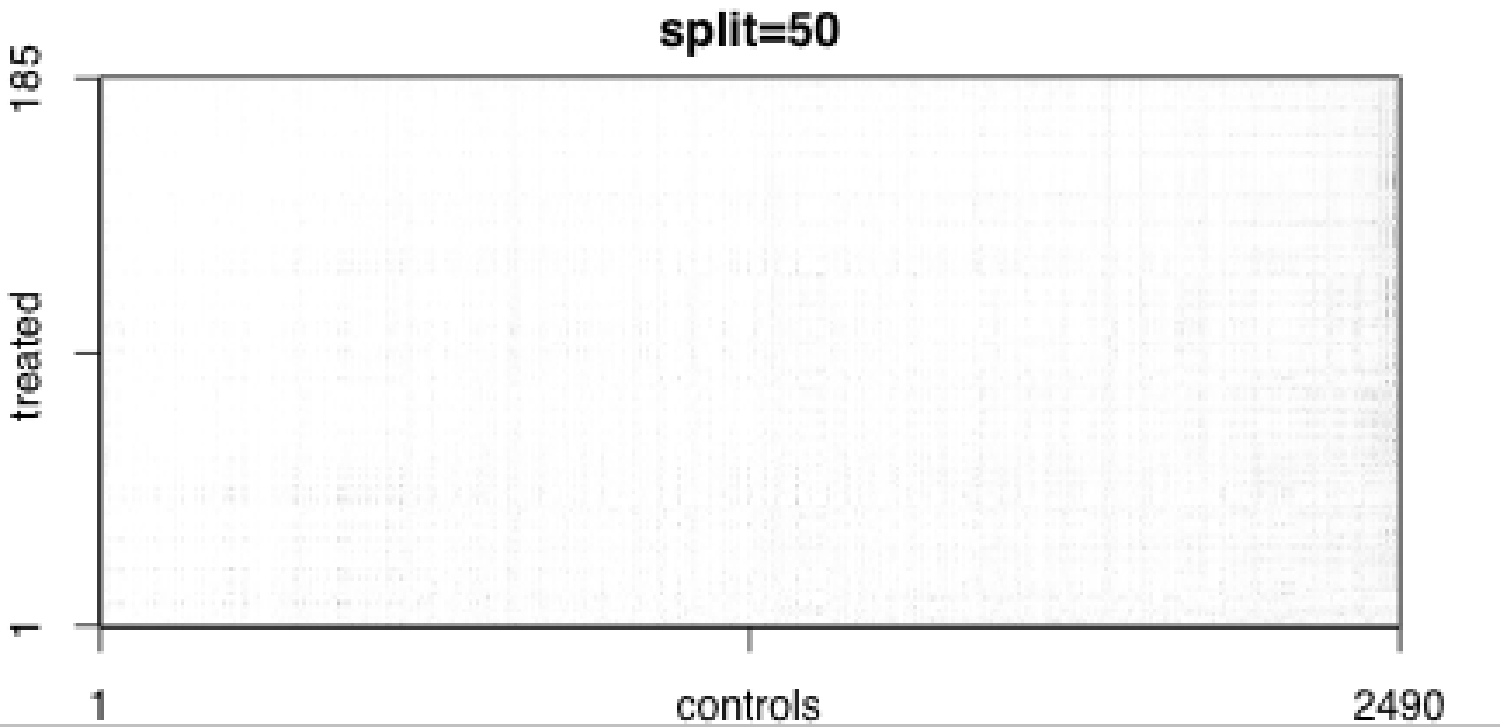}\\
\includegraphics[height=7cm]{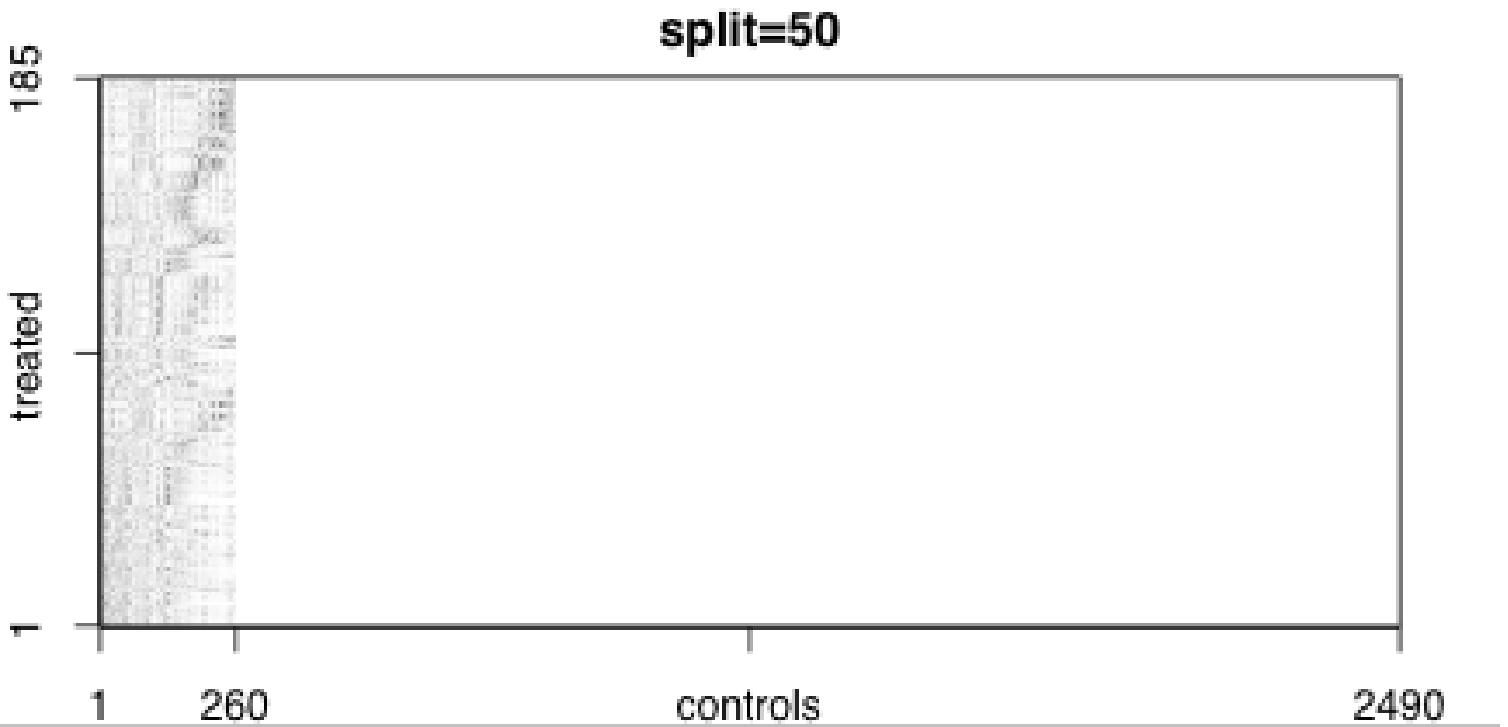}
\end{tabular}
\end{center}
\caption{\small Proximity matrix for random partitions built on the controls. DW experimental dataset (up) and DWvsPSID non experimental dataset (middle). The bottom image is the same as the top one with a 2490-points $x$ axis (as in the middle image). The middle image contains few faint spots, contrary to the top image that shows  the good quality of the match that can be achieved  for the DW experimental data. For further details see \S \ref{sec:proxy}.}
\label{fig2}
\end{figure}

Figure \ref{fig2} reports  a comparison between the proximity matrix for the DW experimental data (upper image) and the DWvsPSID  non experimental data (middle image) with split=50. This value of the split allows for the highest number of matches for the DWvsPSID sample (see Table \ref{tab3}), including spurious ones.
As one can see, the image for DW data is evidence of a good quality matching. On the contrary,
the image for the DWvsPSID data has only few and faint spots, therefore showing that the two groups (treated/controls) in this non experimental sample cannot be successfully matched.
Notice that upper and middle images have different $x$-axis length. One might be tempted to impute to this fact the different intensity in the two images\footnote{In fact, the pixels in the upper image are wider than the ones in the middle image.}. To avoid any misunderstanding, the upper image has been rescaled to the same number of points as the middle image. The result,  reported in the bottom picture of Figure \ref{fig2}, proves that the evidence of the matching quality does not vanish because of rescaling.

This graphical analysis of $\bf P$ is coherent with the analysis based on the indexes used in previous sections to asses the quality of a match.

Once the proximity matrix has been introduced, it seems natural to introduce an ATT estimator based on $\bf P$ defined as follows
\begin{equation}
\tilde\tau =\frac{1}{n_T} \sum_{i\in T}\left(Y_i^T - \frac{1}{\left|C_i\right|}  \sum_{j\in C_i}  p_{ij} Y_j^C \right)
\label{eq:postatt}
\end{equation}
with $C_i = \{ j : p_{ij}>0\}$.
The evaluation of $\tilde\tau$ (and his corrected version $\tilde\tau'$) in our datasets is included in Tables \ref{tab1a}, \ref{tab2a}, \ref{tab3a}. As one can see the results are different (although slightly) from those reported in Tables \ref{tab1}, \ref{tab2} and \ref{tab3}. In fact, the two estimators $\hat\tau$ and $\tilde\tau$ coincide only if
$$
\frac{1}{|C_i|}\sum_{J\in C_i} Y_j^C p_{ij} = 
\frac{1}{R}\sum_{k=1}^R \frac{1}{\left| C_i^{(k)}\right|}\sum_{j\in C_i^{(k)}} Y_j^C w_{ij}^{(k)}
$$
which is true, for instance, with paired sample and split parameter equal to 2.
Also the estimates of the variance are affected by the same difference in the weights.

\section*{Concluding remarks}
The method we propose  to attain a match  directly partitioning the covariates space instead of using a score or a distance approach, seems to be able to discriminate whether two samples can be matched or not, and if the answer is positive, the method provides reliable estimates of the average treatment effect.
Although no formal theory has been developed yet, these computational results are promising.
The idea of representing  the match using the proximity matrix is an additional fruitful tool to understand when a good match cannot be obtained (and hence alternative methods should be used to evaluate the effect of a treatment) and also to see where it failed (as it is clear from the pictures which treated units have been matched and which have not).
The random recursive partitioning method can be used to accomplish several other tasks, like  outliers analysis and clustering classification. Large sample properties of our ATT estimators can be derived as in Abadie and Imbens (2004).

\section*{Acknowledgments}
We wish to thank Robert Gentleman for careful reading of first drafts of the manuscript and for his helpful criticism.

\eject

\begin{landscape}
\begin{table}
\begin{center}
\begin{tabular}{lrrrrr|rrc|rr|cc}
\hline
\hline
&split&$\hat\tau$&$\hat\sigma_{\tau}$ & $\hat\tau'$ & $\hat\sigma_{\tau'}$  & T & C& \% o.t. &SW &SW'&  95\% C.I.($\tau$)&95\% C.I.($\tau'$)\\
\hline
\hline
\multirow{4}{2.5cm}{Lalonde\\297 treated\\ 425 controls\\ naive targets: ATT=886, SD=488}
&50 &822.9 & 492.5 &  958.2 & 644.1 & 297 & 425& -& & $\bullet\bullet$ & (705.0 : 948.1) & (880.0 : 1058.2)\\
&32 &824.1 & 494.2 &  962.3 & 648.2 &  297 & 425& -&&  $\bullet\bullet$ & (682.4 : 970.9) & (859.0 : 1064.8)\\
&20 &820.4 &  495.7  &  956.3  & 653.6  &  296 &  424&-& &  & (643.9 : 993.1) & (805.3 : 1103.0)\\
&16 &820.1 &  498.6  &  955.1 &  659.1 &  296 &  424&  -& & $\bullet$ & (628.1 : 1033.7) & (807.7 : 1132.0)\\
&8 &803.4 &  501.7   & 940.4 &  673.5   & 290  & 424 &  90.8 & & $\bullet\bullet$ & (534.7 : 1097.7) & (728.5 : 1218.3)\\
&4 &755.6 &  463.4  & 907.3 &  658.3   & 267  & 423&  1.2 &$\bullet$ & $\bullet$ &(432.2 : 1150.6) & (598.7 : 1282.2)\\
&2 &875.6 &  253.2   &1070.8 &  535.4  &  235  & 416 &  0.0&  &     & (505.3 : 1299.7) & (714.3 : 1476.3)\\
\hline
\hline
\multirow{4}{2.5cm}{Dehejia-Wahba\\ 185 treated\\ 260 controls\\ naive targets: ATT=1794, SD=670}
&50 &1741.3 & 681.1 &  1718.0 & 680.4  & 185 & 260 & -& &$\bullet\bullet\bullet$ &(1604.0 : 1880.6) & (1601.7 : 1806.0)\\
&32 &1763.1 & 681.7 & 1741.2 & 862.7& 185 & 260 & -& $\bullet$ & $\bullet\bullet$ & (1563.4 : 1917.5) & (1575.8 : 1887.1)\\
&20 &1767.6 &  685.9  & 1753.3  & 827.9 &  185  & 260  & 99.6 &$\bullet\bullet\bullet$& $\bullet$ & (1541.4 : 1979.2) & (1551.1 : 1938.8)\\
&16 &1784.0 &  685.4  & 1777.8 &  876.1  &  184 &  260  &  99.6 & & & (1525.5 : 2052.6) &( 1546.5 : 1991.4)\\
&8 &1750.0 & 691.2 & 1757.6 &  895.4 &  180  & 260 & 88.4 && $\bullet$& (1420.5 : 2124.1) & (1466.1 : 2130.9)  \\
&4 &1691.9 & 641.9 & 1730.6  & 870.4  & 163 &  259 &  1.2 &$\bullet$&  & (1298.2 : 2157.5) & (1336.5 : 2140.1)\\
&2 &1797.3 & 395.9 & 1919.8  & 719.7 & 141 &  254 &   0.0 &  &    & (1288.7 : 2326.1) & (1425.0 : 2473.4)\\
\hline
\hline
\multirow{4}{2.5cm}{Smith-Todd\\108 treated\\ 142 controls\\ naive targets: ATT=2748, SD=1005}
&50& 2651.1 & 1021.1 & 2889.0& 1243.0&  108 &  142 & -& $\bullet\bullet\bullet$&$\bullet\bullet\bullet$& (2481.7 : 2873.0) & (2707.0 : 3019.0)\\
&32& 2657.2 & 1027.0 &   2905.2& 1253.6 &  108 &  142 &-& $\bullet\bullet\bullet$ &$\bullet\bullet\bullet$& (2408.1 : 2876.0 ) & (2716.7 : 3064.2)\\
&20& 2646.6 & 1029.7  & 2905.7 & 1264.4  & 108 &  142 & 98.8 & $\bullet\bullet$& $\bullet\bullet\bullet$& (2340.4 : 2934.4) & (2673.0 : 3138.9)\\
&16& 2648.8 & 1036.1  & 2910.5 & 1277.9 &    108 &  142 &  98.0&  &$\bullet\bullet$ &(2295.7 : 2967.6) & (2610.3 : 3182.7)\\
&8& 2654.9 & 1051.8 &   2912.2 & 1308.3 & 105  & 142 &  79.2 &$\bullet\bullet$ &  & (2096.4 : 3273.9) & (2410.3 : 3413.7)\\
&4& 2705.0 &  976.1 & 2945.5 &  1261.8  &  96  & 141 & 10.4 &$\bullet\bullet$& $\bullet\bullet$ & (1997.9 : 3524.4) & (2401.9 : 3625.2)\\
&2& 2680.2 &  319.8 & 3014.0 &  859.8&  82 &  137 &   0.0&  &   & (1858.0 : 3499.3) & (2181.6 : 3917.3)\\
\hline
\end{tabular}
\caption{\small The results of random trees  built on the treated. Average values over 250 replications. $\hat\tau$ and $\hat\sigma_{\tau}$ are the estimators of $\tau$ and its standard deviation in the population ($\hat\tau'$ and  $\hat\sigma_{\tau'}$ are the difference-in-covariate corrected versions, see \S \ref{sec:idea}). $T$ and $C$ are respectively the average number of treated and controls matched per tree. The percentage of trees (``-'' = 100\%) that match at least 95\% of the treated is reported in column ``\% o.t.''.
SW and SW' report the results of the Shapiro-Wilk test for normality respectively for $\hat\tau$ and $\hat\tau'$: the bullets ($\bullet$) mean that the hypothesis of normality is rejected at the corresponding level ($\bullet=0.05$, $\bullet\bullet=0.01$, $\bullet\bullet\bullet=0.001$). The last two columns are the 95\% Monte Carlo confidence intervals for $\tau$ and $\tau'$.}
\label{tab1}
\end{center}
\end{table}
\end{landscape}

\begin{landscape}
\begin{table}
\begin{center}
\begin{tabular}{lrrrrr|rrc|rr|cc}
\hline
\hline
&split&$\hat\tau$&$\hat\sigma_{\tau}$ & $\hat\tau'$ & $\hat\sigma_{\tau'}$  & T & C&\% o.t. &SW &SW' &  95\% C.I.($\tau$)&95\% C.I.($\tau'$)\\
\hline
\hline
\multirow{4}{2.5cm}{Lalonde\\297 treated\\ 425 controls\\ naive targets: ATT=886, SD=488}
 & 50&  820.7&   493.2&  944.0&    645.5&  297 &  424&  -& $\bullet$ &$\bullet$& (692.7 : 955.7) & (834.6 : 1062.0)\\
 & 32& 823.6&    494.5&    940.9&   648.7&  297 &  424& -&     &   & (667.0 : 956.4) & (811.7 : 1051.7)\\
& 20&  824.8&   498.2&  935.8&   655.9&  297 &  421 & -&$\bullet$ &   & (667.6 : 992.6) & (790.1 : 1096.6)\\
& 16&  836.7&    499.1&     949.2&     659.7 &  297 &  417 &-&       &   & (638.0 : 1022.1 ) & (763.4 : 1114.7)\\
& 8& 852.6&    501.9&     958.2&    673.8&  297 &  388&  - &  &  & (588.0 : 1145.3) & (702.9 : 1203.5)\\
& 4&  849.8&     475.6&    950.1&   679.1&   297 &  322&  -&  &    & (467.5 : 1230.9) & (570.7 : 1316.7)\\
& 2&  790.5&    306.8&   849.0 &   610.5 &  295 &  247&  -&  &   & (292.9 : 1299.4) & (375.8 : 1344.4)\\
\hline
\hline
\multirow{4}{2.5cm}{Dehejia-Wahba\\ 185 treated\\ 260 controls\\ naive targets: ATT=1794, SD=670}
& 50& 1733.3&  682.9&  1701.4&   860.3&  185 &  260& -&   $\bullet\bullet\bullet$& $\bullet\bullet\bullet$& 1571.7 : 1867.7) & (1558.4 : 1827.8)\\
& 32& 1759.1&   684.2&  1730.9&   865.5 &  185 &  259 &-& $\bullet\bullet\bullet$& $\bullet\bullet\bullet$ & 1580.7 : 1939.4) & (1568.3 : 1897.1)\\
& 20&  1773.3&  688.6&   1753.0 &  873.4&  185 &  259& -&  $\bullet$&  & 1531.7 : 1968.6) & (1539.7 : 1938.6)\\
& 16&  1785.8&  688.7&   1761.0&  877.5&  185 &  256&  -&   &  & (1562.2 : 2015.3) & (1567.6 : 1954.7)\\
& 8& 1808.4&   692.4&   1776.1&  895.3& 185& 239&  -&  &    & (1525.5 : 2123.5) & (1535.8 : 2059.3)\\
& 4& 1826.7&   651.3&   1772.4&  891.1&  185 &  200 & -&       &  & (1383.7 : 2267.6) & (1399.6 : 2187.9)\\
& 2&  1770.2&  474.7&  1620.2& 801.8 &  183 &  164 &  -&   & $\bullet$& (1221.0 : 2407.0) & (1171.7 : 2228.4)\\
\hline
\hline
\multirow{4}{2.5cm}{Smith-Todd\\108 treated\\ 142 controls\\ naive targets: ATT=2748, SD=1005}
& 50 & 2636.6 & 1020.5&  2869.1 & 1243.1&  108 &  142& -& $\bullet\bullet$ &$\bullet$  & (2467.1 : 2855.9) & (2740.4 : 2995.1)\\
& 32&  2622.6&   1028.7& 2857.7&   1255.8 & 108 &  142 &  -&  & $\bullet\bullet$ & (2386.7 : 2914.1) & (2629.3 : 3062.1)\\
& 20&  2859.4&   1044.3 &  2834.0 & 1277.5 & 108 &  140 &  -&  &     & (2279.3 : 2884.0) & (2570.7 : 3087.8)\\
& 16 &  2622.6&   1043.0&    2850.2&   1281.5&  108 &  138 & -&   & $\bullet$& (2265.7 : 2966.1) & (2570.3 : 3202.6)\\
& 8 & 2859.3&   1055.2&   2836.6&  1315.6&  108 &  125 &   -&   &   & (2108.9 : 3155.9) & (2355.0 : 3310.1)\\
& 4&  2620.6&    1001.3&  2756.9 &   1307.7&  108 &  101 &  -&    &  & (1801.1 : 2258.4) & (2148.7 : 3412.9)\\
& 2&  2709.2&    598.1&  2691.9 & 1079.2&  107 &  75&  -&  $\bullet\bullet$ &  & (1809.2 : 3724.0) & (1862.1 : 3548.2)\\
\hline
\end{tabular}
\caption{\small The results of random trees  built on the controls. Average values over 250 replications. $\hat\tau$ and $\hat\sigma_{\tau}$ are the estimators of $\tau$ and its standard deviation in the population ($\hat\tau'$ and  $\hat\sigma_{\tau'}$ are the difference-in-covariate corrected versions, see \S \ref{sec:idea}). $T$ and $C$ are respectively the average number of treated and controls matched per tree. The percentage of trees (``-'' = 100\%) that match at least 95\% of the treated is reported in column ``\% o.t.''.
SW and SW' report the results of the Shapiro-Wilk test for normality respectively for $\hat\tau$ and $\hat\tau'$: the bullets ($\bullet$) mean that the hypothesis of normality is rejected at the corresponding level ($\bullet=0.05$, $\bullet\bullet=0.01$, $\bullet\bullet\bullet=0.001$). The last two columns are the 95\% Monte Carlo confidence intervals for $\tau$ and $\tau'$.}
\label{tab2}
\end{center}
\end{table}
\end{landscape}

\begin{landscape}
\begin{table}
\begin{center}
\begin{tabular}{lrrrr|rrc|rr|cc}
\hline
\hline
split&$\hat\tau$&$\hat\sigma_{\tau}$ & $\hat\tau'$ & $\hat\sigma_{\tau'}$ & T & C& \% o.t. &SW &SW'&  95\% C.I.($\tau$)&95\% C.I.($\tau'$)\\
\hline
50 & -1931.1 & 1357.4 & 557.2 & 1848.5 & 175 & 2398 & 54.0 &   $\bullet\bullet\bullet$ &  $\bullet\bullet\bullet$  & (-3645.7 :  312.1) & (-669.4 :  2432.8)\\
32 & -1663.9 & 1198.2 & 806.4 & 1697.6 & 169 & 2359 & 26.0 & $\bullet\bullet\bullet$ &  $\bullet\bullet\bullet$ & (-3200.9 :  254.2) & (-296.3 :  2191.2)\\
20 & -1749.8 & 1117.2 & 752.7 & 1610.6 & 161 & 2323 & 2.4 & $\bullet\bullet$ &  $\bullet\bullet$  & (-3465.2 :  390.9) & (-463.5 :  2314.9)\\
16 & -1777.2 & 1058.0 & 789.2 & 1569.5 & 159 & 2306 & 1.2 & $\bullet\bullet$   &   & (-3354.4 : -71.2) & (-414.4 : 1997.4)\\
8 & -2127-1 & 938.1 & 755.5 & 1439.9 & 148 & 2221 & 0.0 &    &      & (-3884.0 : -351.4) & (-302.5 : 1848.2)\\
4 & -2604.8 & 761.3 & 813.5 & 1280.9 & 146 & 1915 &  0.0 &  &  $\bullet$  & (-4487.6 :  -780.0) & (-270.0 : 2138.5)\\
2 & -2736.9 & 491.3 & 612.2 & 1139.5 & 146 & 349 & 0.0 &     &   & (-3940.6 : -1269.3) & (-417.0 : 1697.7)\\
\hline
\\
\\
\hline
  50&  -692.4 &   1077.0&  1049.7&    1657.8&  176 &  552&  58.4& $\bullet\bullet$   &$\bullet\bullet\bullet$ & (-2271.0 :  1255.0) & (5.1 : 2318.1)\\
32 & -619.3 & 1072.2 & 984.4 & 1636.7 & 177 & 449 & 58.4 & $\bullet\bullet\bullet$   & $\bullet\bullet\bullet$  & (-2364.2 :  1222.4) & (-208.2 :  2401.3)\\
20 & -392.9 & 1034.6 & 1058.3 & 1588.2 & 177 & 359 & 54.4 &$\bullet\bullet\bullet$ & $\bullet\bullet\bullet$  & (-1947,2 :  1425.2) & (-39.9 :  2472.2) \\
16 & -295.2 & 1033.0 & 1065.1 & 1586.1 & 177 & 321 & 55.2 & $\bullet\bullet\bullet$   & $\bullet\bullet\bullet$  & (-1583.8 : 1398.2) &  (34.9 :  2466.6)\\
8 & -64.7  & 977.7 & 1263.8 & 1572.6 & 177 & 243 & 51.2 & $\bullet\bullet$  & $\bullet\bullet\bullet$  & (-1481.8 : 1373.7) &  (221.5 : 2499.6) \\
4 & -877.6 & 997.7 & 1211.7 & 1643.1 & 176 & 178 & 52.4 & $\bullet\bullet$   & $\bullet$ & (-2610.1 : 741.2) & (80.9 : 2572.7)\\
2 & -323.4 & 534.5 & 1287.6 & 1452.1 & 166 & 128 & 6.4 & $\bullet\bullet$   &  & (-1663.9 : 1101.4) & (180.0 :  2494.1)\\
\hline
\end{tabular}
\caption{\small Random tree results from the tree built on the treated (up) and on the controls (down) for the DW {\it versus} PSID sample.  Average values over 250 replications. $\hat\tau$ and $\hat\sigma_{\tau}$ are the estimators of $\tau$ and its standard deviation in the population ($\hat\tau'$ and  $\hat\sigma_{\tau'}$ are the difference-in-covariate corrected versions, see \S \ref{sec:idea}). $T$ and $C$ are respectively the average number of treated and controls matched per tree. The percentage of trees   that match at least 95\% of the treated is reported in column ``\% o.t.''.
SW and SW' report the results of the Shapiro-Wilk test for normality respectively for $\hat\tau$ and $\hat\tau'$: the bullets ($\bullet$) mean that the hypothesis of normality is rejected at the corresponding level ($\bullet=0.05$, $\bullet\bullet=0.01$, $\bullet\bullet\bullet=0.001$). The last two columns are the 95\% Monte Carlo confidence intervals for $\tau$ and $\tau'$.}
\label{tab3}
\end{center}
\end{table}
\end{landscape}

\begin{table}
\begin{center}
\begin{tabular}{lrrrrr|rrc}
\hline
\hline
&split&$\tilde\tau$&$\tilde\sigma_{\tau}$ & $\tilde\tau'$ & $\tilde\sigma_{\tau'}$  & T & C& \% o.t. \\
\hline
\hline
\multirow{4}{2.5cm}{Lalonde\\297 treated\\ 425 controls\\ naive targets: ATT=886, SD=488}
&50 &866.0 & 402.8 &  974.7 & 423.5 & 297 & 425& -\\
&32 &849.5 & 403.2 &  965.5 & 423.9 &  297 & 425& -\\
&20 &848.8 &  403.8  &  966.1  & 424.6  &  296 &  424&-\\
&16 &851.1 &  404.2  &  966.1 &  425.1 &  296 &  424&  -\\
&8 &841.9 &  406.4   & 962.5 &  427.4   & 290  & 424 &  90.8 \\
&4 &829.4 &  409.7  & 956.3 &  431.1   & 267  & 423&  1.2 \\
&2 &851.1 &  420.5   &1006.0 &  422.9  &  235  & 416 &  0.0\\
\hline
\hline
\multirow{4}{2.5cm}{Dehejia-Wahba\\ 185 treated\\ 260 controls\\ naive targets: ATT=1794, SD=670}
&50 &1766.2 & 579.7 &  1794.0 & 599.5  & 185 & 260 & -\\
&32 &1759.6 & 580.1 & 1708.7 & 600.0& 185 & 260 & -\\
&20 &1737.7 &  581.0  & 1708.6  & 600.9 &  185  & 260  & 99.6\\
&16 &1742.0 &  581.5  & 1727.2 &  601.5  &  184 &  260  &  99.6 \\
&8 &1725.8 & 584.5 & 1730.6 &  604.7 &  180  & 260 & 88.4 \\
&4 &1685.3 & 589.6 & 1685.4  & 610.3  & 163 &  259 &  1.2 \\
&2 &1720.8 & 603.4 & 1756.5  & 625.2 & 141 &  254 &   0.0 \\
\hline
\hline
\multirow{4}{2.5cm}{Smith-Todd\\108 treated\\ 142 controls\\ naive targets: ATT=2748, SD=1005}
&50& 2718.2 & 870.7 & 2873.3& 914.4&  108 &  142 & -\\
&32& 2692.8 & 871.4 &   2876.4& 915.2 &  108 &  142 &-\\
&20& 2655.3 & 872.3  & 2875.4 & 916.3  & 108 &  142 & 98.8\\
&16& 2662.7 & 872.9  & 2879.0 & 917.0 &    108 &  142 &  98.0\\
&8& 2627.0 & 875.7 &   2884.6 & 920.3 & 105  & 142 &  79.2 \\
&4& 2620.7 &  879.2 & 2882.4 &  924.4  &  96  & 141 & 10.4 \\
&2& 2560.8 &  898.1 & 2881.1 &  946.6&  82 &  137 &   0.0\\
\hline
\end{tabular}
\caption{\small The results of random trees  built on the treated over 250 replications. $\tilde\tau$ and $\tilde\sigma_{\tau}$ are the estimators of $\tau$ and its standard deviation in the population ($\tilde\tau'$ and  $\tilde\sigma_{\tau'}$ are the difference-in-covariate corrected versions, see \S \ref{sec:idea}). $T$ and $C$ are respectively the average number of treated and controls matched per tree. The percentage of trees (``-'' = 100\%) that match at least 95\% of the treated is reported in column ``\% o.t.''.}
\label{tab1a}
\end{center}
\end{table}

\begin{table}
\begin{center}
\begin{tabular}{lrrrrr|rrc}
\hline
\hline
&split&$\tilde\tau$&$\tilde\sigma_{\tau}$ & $\tilde\tau'$ & $\tilde\sigma_{\tau'}$  & T & C&\% o.t. \\
\hline
\hline
\multirow{4}{2.5cm}{Lalonde\\297 treated\\ 425 controls\\ naive targets: ATT=886, SD=488}
 & 50&  842.5&   403.1&  947.5&    423.8&  297 &  424&  -\\
 & 32& 847.3&    403.7&    949.4&   424.5&  297 &  424& -\\
 & 20& 836.3  & 404.6   &  943.3  & 425.5   &  297 &  421 & -\\
 & 16&   838.2&   405.2  &   949.5   &   426.2    &  297 &  417 &-\\
 & 8&  856.8&  408.4   & 962.4     &  429.6   &  297 &  388&  -\\
 & 4&   862.8&   413.4   &  961.5   &   435.2 &   297 &  322&  -\\
 & 2&  892.6 &    428.9 &  903.8   &  452.1   &  295 &  247&  -\\
\hline
\hline
\multirow{4}{2.5cm}{Dehejia-Wahba\\ 185 treated\\ 260 controls\\ naive targets: ATT=1794, SD=670}
& 50& 1753.6&  580.1&  1702.3&   600.0&  185 &  260& -\\
& 32& 1758.9&   580.8&  1708.0&   578.4 &  185 &  259 &-\\
& 20&  1764.8&  582.0&   1737.0 &  602.0&  185 &  259& -\\
& 16&  1790.6&  582.8&   1758.6&  602.9&  185 &  256&  -\\
& 8& 1801.0&   587.0&   1760.0&  607.5& 185& 239&  -\\
& 4& 1818.2&   595.0&   1781.9&  616.1&  185 &  200 & -\\
& 2&  1845.6&  613.7&  1644.9& 636.4 &  183 &  164 &  -\\
\hline
\hline
\multirow{4}{2.5cm}{Smith-Todd\\108 treated\\ 142 controls\\ naive targets: ATT=2748, SD=1005}
& 50 & 2671.5 & 871.1&  2862.3 & 914.9&  108 &  142& -\\
& 32&  2627.3&   871.8& 2830.1&   915.8 & 108 &  142 &  -\\
& 20&  2603.4&   872.1 &  2825.4 & 917.3 & 108 &  140 &  -\\
& 16 &  2628.7&   874.0&    2851.2&   918.3&  108 &  138 & -\\
& 8 & 2610.5&   877.5&   2846.7&  922.5&  108 &  125 &   -\\
& 4&  2589.7&    883.4&  2750.4 &   929.4&  108 &  101 &  -\\
& 2&  2883.2&    913.2&  2787.4 & 912.0&  107 &  75&  -\\
\hline
\end{tabular}
\caption{\small The results of random trees  built on the controls over 250 replications.  $\tilde\tau$ and $\tilde\sigma_{\tau}$ are the estimators of $\tau$ and its standard deviation in the population ($\tilde\tau'$ and  $\tilde\sigma_{\tau'}$ are the difference-in-covariate corrected versions, see \S \ref{sec:idea}). $T$ and $C$ are respectively the average number of treated and controls matched per tree. The percentage of trees (``-'' = 100\%) that match at least 95\% of the treated is reported in column ``\% o.t.''.}
\label{tab2a}
\end{center}
\end{table}

\begin{table}
\begin{center}
\begin{tabular}{lrrrr|rrc}
split&$\tilde\tau$&$\tilde\sigma_{\tau}$ & $\tilde\tau'$ & $\tilde\sigma_{\tau'}$ & T & C& \% o.t. \\
\hline
50 & -3592.2 & 597.6 & 747.6 & 665.5 & 175 & 2398 & 54.0 \\
32 & -3381.2 & 602.1 & 980.5 & 671.8 & 169 & 2359 & 26.0 \\
20 & -3599.8 & 606.9 & 937.7 & 678.4 & 161 & 2323 & 2.4\\
16 & -3756.8 & 612.6 & 1100.8 & 686.6 & 159 & 2306 & 1.2\\
8 & -4441.8 & 627.3 & 1126.8 & 706.5 & 148 & 2221 & 0.0 \\
4 & -6516.3 & 620.7 & 1066.1 & 697.5 & 146 & 1915 &  0.0\\
2 & -3380.3 & 658.7 & 906.9 & 11739.9 & 146 & 349 & 0.0 \\
\\
\\
split&$\tilde\tau$&$\tilde\sigma_{\tau}$ & $\tilde\tau'$ & $\tilde\sigma_{\tau'}$ & T & C& \% o.t. \\
\hline
  50&  -649.3& 594.7 &1162.8 & 659.8  &  176 &  552&  58.4\\
 32& -428.2 & 593.3 &1131.3 &657.3   & 177 & 449 & 58.4 \\
 20&-137.7  &  596.0&1162.7 &  661.3 &177 & 359 & 54.4 \\
 16 &  -5.2&  598.5& 1210.8&   664.8& 177 & 321 & 55.2 \\
 8 & -42.7 &618.9  & 1298.2& 693.5  &177 & 243 & 51.2 \\
 4 & -1457.0 & 643.0 &1263.0 &728.4   &176 & 178 & 52.4\\
 2 &-1167.5  & 731.8 & 1326.4& 842.8  & 166 & 128 & 6.4\\
\hline
\end{tabular}
\caption{\small Random tree results over 250 replications from the tree built on the treated (up) and on the controls (down) for the DW {\it versus} PSID sample. $\tilde\tau$ and $\tilde\sigma_{\tau}$ are the estimators of $\tau$ and its standard deviation in the population ($\tilde\tau'$ and  $\tilde\sigma_{\tau'}$ are the difference-in-covariate corrected versions, see \S \ref{sec:idea}). $T$ and $C$ are respectively the average number of treated and controls matched per tree. The percentage of trees   that match at least 95\% of the treated is reported in column ``\% o.t.''.}
\label{tab3a}
\end{center}
\end{table}

\end{document}